\documentclass[11pt]{amsart}
\usepackage{latexsym}
         \usepackage[all]{xy}
                \usepackage{amscd}
                \usepackage{amssymb}
                \usepackage{amsmath}
                \usepackage{amsthm}
                \usepackage{enumerate}
                \usepackage{verbatim}

\def\co{\colon\thinspace}
\def\d{\delta}
\def\e{\epsilon}

\def\G{\Gamma}
\def\b{\beta}

\def\ga{\gamma}

\def\a{\alpha}%
\def\"#1{{\accent"7F #1\penalty10000\hskip 0pt plus 0pt}} 

\def\N{{\Bbb N}}
\def\Z{{\Bbb Z}}
\def\R{{\Bbb R}}

\def\sp{{\frak {sp}}}

\def\h{{\frak h}}

\def\m{{\frak m}}

\newtheorem{thm}{Theorem}[section]

\newtheorem{quest}[thm]{Question}
\newtheorem{lem}[thm]{Lemma}
\newtheorem{prop}[thm]{Proposition}

\newtheorem{Example}[thm]{Example}

\newtheorem{Counterexample}[thm]{Counterexample}

\newtheorem{remark}[thm]{Remark}
\newenvironment{rmk}{\begin{remark}\rm}{\end{remark}}
\newtheorem{Fact}[thm]{Fact}

\newtheorem{Nothing}[thm]{$\!\!\!$}

\newcommand{\be}{\begin{equation} }

\newcommand{\ene}{\end{equation} }

\newcommand{\ba}{\begin{eqnarray}}
\newcommand{\ea}{\end{eqnarray}}
\newcommand{\ban}{\begin{eqnarray*}}
\newcommand{\ean}{\end{eqnarray*}}

\newcommand{\Aut}{\mbox{\rm Aut}}

\newcommand{\Id}{\mbox{\rm Id}}
\newcommand{\Diff}{\mbox{\rm Diff}}
\newcommand{\diam}{\mbox{\rm diam}}
\newcommand{\Ric}{\mbox{\rm Ric}}

\newtheorem{main}{Theorem}

\begin{document}
\abovedisplayskip=6pt plus3pt minus3pt \belowdisplayskip=6pt
plus3pt minus3pt
\title[Nonnegative pinching, moduli spaces and bundles with infinitely 
many souls]{\bf Nonnegative pinching, moduli spaces and \\
 bundles with infinitely many souls}

\rm
\thanks{\it The first author was supported in part by the NSF grant \# DMS-
0204187, the last author by a DFG Heisenberg Fellowship}
\author{Vitali Kapovitch}
\address{Vitali Kapovitch\\Department of Mathematics\\University of M\"unster\\
48149 M\"unster, Germany}\email{kapowitc@math.uni-muenster.de}
\author{ Anton Petrunin }\address{Anton Petrunin\\ Department of Mathematics\\ Pennsylvania State University\\
University Park, State College, PA 16802
}\email{petrunin@math.psu.edu}
\author{Wilderich Tuschmann}\address{Wilderich Tuschmann
\\Department of Mathematics\\University of M\"unster\\
48149 M\"unster, Germany}\email{wtusch@math.uni-muenster.de}

\maketitle
\begin{abstract}
We show that in each dimension $n\ge 10$ there
exist infinite sequences of homotopy equivalent
but mutually non-homeomorphic closed simply connected 
Riemannian $n$-manifolds with $0\le \sec\le 1$, positive Ricci curvature 
and uniformly bounded diameter.
We also construct open manifolds of fixed diffeomorphism type
which admit infinitely 
many complete nonnegatively pinched metrics 
with souls of bounded diameter 
such that the souls are mutually non-homeomorphic. 
Finally, we construct examples of noncompact manifolds
whose moduli spaces of complete metrics with $\sec\ge 0$ 
have infinitely many connected components.
\end{abstract}

\section{Introduction}


\noindent

In this 
article we discuss several infiniteness phenomena 
in 
nonnegative sectional curvature.

Our first such result is motivated by the finiteness theorems in
Riemannian geometry and a question of S.-T.~Yau
which asks
whether there always exists   only a finite number of
diffeomorphism types of closed smooth manifolds of positive
sectional curvature that are homotopy equivalent to a given
positively curved manifold (\cite{Yau}, Problem~11).

If one relaxes the condition $\sec>0$ to $\sec\ge 0$ 
then the answer to Yau's question is known to be false in all
dimensions $\ge 7$,
even in the category of simply connected manifolds.
Counterexamples can  here be obtained 
in the following way: By 
a result of Grove and Ziller~\cite{GZ}
the total space of any linear $S^k$-bundle over $S^4$ admits 
a Riemannian metric with nonnegative curvature.
However,  for $k\ge 3$
the total spaces of such bundles  fall into infinitely many
homeomorphism, but only finitely many homotopy types 
(\cite{DW}; if $k=3$, one has in addition to assume that the
Euler class of the bundle be zero).
On the other hand, 
when rescaled to have uniformly bounded diameter,
by \cite{PT} these  examples 
cannot satisfy any uniform upper curvature bound. 
More generally,  it is natural to look at the following question:

\begin{quest}\label{intro-quest1}
Given fixed $n\in \N, D>0$ and $c, \, C\in \R$,
are there at most finitely many diffeomorphism classes 
of pairwise homotopy equivalent closed Riemannian 
$n$-manifolds $M^n$
with sectional curvature $c\le\sec\le C$ and diameter $\le D$?
\end{quest}

The diffeomorphism finiteness theorems in Riemannian geometry
(see, e.g.,  \cite{AC,Ch1,Pet,FR1,FR3,GPW,PT,Tu})
leave this question in general dimensions completely open.
However, the answer is known to be positive in some special situations. 

This is, for example, the case
if $c>0$ and $n=2m$ by \cite{Kl, Ch1},
if $M^n$, $n\ne 4$, is simply connected and $C=4c>0$ by \cite{Ber, CG},
if $D=D(C,c,n)$ is sufficiently small 
by Gromov's theorem on almost flat manifolds \cite{Gr,BuKa,Ruh}  
and the rigidity of infranilmanifolds~\cite{Au} (cf. ~\cite{FH}),
if $M$ is $2$-connected by ~\cite{PT},  
or if $C\le 0$ and $n\ge 5$ by results 
of Farrell and Jones~\cite{FJ1,FJ2}. 
Remarkably enough, 
in the latter case one actually does not even need the lower curvature 
and the  upper diameter bounds. In other words, 
for $n\ge 5$ the answer to the analogue of Yau's question 
for nonpositive curvature 
(which in this case is a special case of the Borel conjecture) is yes.

As a preliminary result we first show that in general  the answer to Question~\ref{intro-quest1} 
is  actually negative in all dimensions $\ge 7$:

\begin{prop}\label{intro-curv} There exists $D>0$  such that
for any $n\ge 7$ there exist  an infinite sequence of  homotopy equivalent 
but mutually non-homeomorphic closed Riemannian $n$-manifolds $M^n_k$ with
$$| \sec  (M^n_k) |\le 1 \,
\text{ and } \, \diam (M^n_k) \le D.$$
If $n\ne 8$, all these manifolds can in addition be chosen to be simply-connected.
\end{prop}

Notice that for simply connected manifolds,  by \cite{FR1, Tu},   
$n=7$ is indeed the smallest dimension where such sequences can occur.

Our first main concern in this paper is, however,
the analogue of  Yau's question for nonnegative pinching, i.e.,
the following special case of Question~\ref{intro-quest1}:

\begin{quest}\label{intro-quest2}
Given fixed $n\in \N$ and $C,D>0$, 
are there always at most finitely many diffeomorphism types 
of pairwise homotopy equivalent closed Riemannian $n$-manifolds 
with sectional curvature $0\le\sec\le C$ and diameter $\le D$?

\end{quest}



Notice here that starting from dimension $n=6$, 
from~\cite{GZ} one may infer the existence
of infinite sequences of 
closed simply connected nonnegatively curved $n$-manifolds
of mutually distinct homotopy type,
and in dimensions $n>8$, $n\ne 10$
by \cite{To} 
there even exist infinite sequences of 
closed simply connected nonnegatively pinched Riemannian $n$-manifolds
with pairwise non-isomorphic rational cohomology rings
that also satisfy uniform upper diameter bounds.
Totaro (\cite{To}) also showed
that there exist infinite sequences of 
closed simply connected nonnegatively curved $6$-manifolds
with pairwise non-isomorphic rational cohomology rings,
and (improving earlier work of \cite{FR2} for manifolds of dimension $n\ge 22$)
that for any $n\ge 7, n\ne 8$ there exists   infinitely many  
closed simply connected Riemannian $n$-manifolds
with $|sec|\le 1$ and uniformly bounded diameter and
pairwise non-isomorphic rational cohomology rings.

Our first main result
shows that if $n\ge 10$, the answer to Question~\ref{intro-quest2} 
is in general negative, even under the extra assumption
of  positive Ricci curvature.

\begin{main}\label{intro-main}
There  exists  $D>0$  such that for each dimension $n\ge10$ there exists an  infinite sequence
$(M^n_k)_{k\in\N}$ of 
pairwise homotopy equivalent but mutually non-homeomorphic 
closed simply connected 
Riemanni\-an $n$-manifolds satisfying

$$0\le \sec  (M^n_k) \le 1, \,\Ric(M^n_k) >0 \,
\text{ and } \, \diam (M^n_k) \le D.$$

\end{main}

Notice that this result relates Yau's problem to  another 
long standing open question in Riemannian geometry:

Are there any obstructions to the existence of a Riemannian metric 
with positive sectional curvature on a closed  simply connected manifold
of nonnegative sectional and positive Ricci curvature?


It is quite likely that the dimensional restriction $n\ge10$ in Theorem $A$
is not optimal.
Again, by \cite{FR1,Tu}, this dimension must be at least $7$.

We continue with a description of the further main results
of this paper, which concern
the souls and moduli spaces of metrics of
open manifolds of nonnegative sectional curvature.

Very recently, in~\cite{Bel} Belegradek constructed 
the first examples of manifolds 
admitting infinitely many nonnegatively curved metrics with 
mutually non-homeomorphic souls.
In our second main theorem  we  sharpen this result 
by constructing such examples 
which in addition have uniform  bounds on the curvature 
of the manifolds and the diameters of the souls:

\begin{main}\label{intro-main2}
For any $k>10$ 
the manifold $S^2\times S^2\times S^3\times S^3\times \R^k$
admits an infinite sequence of complete nonnegatively curved 
metrics $g_i$  with pairwise non-homeomorphic souls $S_i$ such that
\[
  0\le \sec(M,g_i)\le 1 \,
\text{ and } \,  \diam (S_i)\le D
\]
where $D$ is a positive constant independent of $k$ and $i$.
\end{main}

Another interesting application of the construction that we employ
in the proof of Theorem~\ref{intro-main2}
concerns the moduli spaces of open manifolds with 
nonnegative sectional curvature.
Let $\mathfrak R_{\sec \ge 0}(M)$
denote the space of complete Riemannian metrics
with nonnegative sectional curvature on a given smooth manifold $M$.
Then Diff(M), the group of diffeomorphisms of $M$,
acts on this space by pulling back metrics,
and the orbit space $\mathfrak R_{\sec \ge 0}(M)/\Diff(M)$
is called the moduli space of (complete) nonnegatively
curved Riemannian metrics on $M$.
We show:

\begin{main}\label{intro-main3}
There exist a manifold $M^{22}$  which admits an infinite 
sequence of complete metrics $g_i$  with pairwise non-homeomorphic souls $S_i$ such that
\[
  0\le \sec(M,g_i)\le 1 \,
\text{ and } \,  \diam (S_i)\le D
\]

and such that the equivalence classes of the metrics $g_i$ 
all lie in different connected components
of the moduli space $\mathfrak R_{\sec \ge 0}(M)/\Diff(M)$ 
of complete metrics with $\sec\ge 0$ on $M$.

Moreover, for any closed nonnegatively curved manifold $(N,g)$,
the product metrics $g_i\times g$ 
all lie in different connected components of  
the moduli space $\mathfrak R_{\sec \ge 0}(M\times N)/\Diff(M\times N)$.
\end{main}

To put Theorem~\ref{intro-main3} into further perspective,
we note first
that in~\cite{KrSt}  Kreck and Stoltz constructed  a closed manifold $M^7$
such that the moduli space of metrics of positive Ricci curvature on $M$ 
has infinitely many connected components.  In fact,  by choosing somewhat 
different metrics their methods also show that  
$\mathfrak R_{\sec \ge 0}(M)/\Diff (M)$ 
also has infinitely many connected components!  
Since this was not observed in~\cite{KrSt}, 
let us briefly explain why that is true.

Kreck and Stoltz actually construct an invariant $s$ which 
distinguishes components of  $\mathfrak R_{\rm{scalar} >0}(M)/\Diff(M)$. 
They look at various $S^1$-bundles  over $S^2\times CP^2$ with indivisible Euler classes, which admit Einstein metrics of positive Ricci curvature 
constructed by Wang and Ziller(~\cite{WZ}).
It is then  shown~\cite[Theorem 3.11]{KrSt} 
that provided the metrics are $S^1$ invariant and have totally geodesic 
fibers (which is true for Wang-Ziller metrics), 
the invariant $s$ depends only on the Euler class of the bundle.  
One can exhibit infinitely many bundles with distinct $s$ 
invariants but diffeomorphic total spaces~\cite[Theorem 3.2, 3.4]{KrSt}.  
Unfortunately, the Einstein metrics given by ~\cite{WZ} 
do not have nonnegative sectional curvature. 
However, we notice here that  one can represent 
any $S^1$ bundle over $S^2\times CP^2$ with an indivisible Euler class 
as a free isometric quotient  $(S^3\times S^5)/ S^1$. 
The natural Riemannian submersion metric coming from the product metric
 on $S^3\times S^5$ is easily seen to have $\sec\ge 0$ and $\Ric>0$, 
it has totally geodesic fibers and it is $S^1$-invariant. 
Therefore the same bundles as considered in \cite{KrSt},
but taken with these metrics have distinct $s$-invariants 
and hence lie in different components of 
$\mathfrak R_{\rm{scalar} >0}(M)/\Diff(M)$.  Since any metric of $\sec\ge 0$ on $M$ has scalar$\ge 0$ and scalar $>0$ at a point, by~\cite{Aub} it can be deformed to a nearby metric of scalar$>0$. Therefore the above metrics  also lie in different components of $\mathfrak R_{\sec\ge 0}(M)/\Diff(M)$.

Observe, however,  
that all the different components of 
$\mathfrak R_{\rm{scalar}> 0}(M)/\Diff(M)$ obviously become connected 
if we stabilize $M$ by multiplying it by a closed manifold 
with nonnegative sectional and positive scalar curvature, 
for example by $S^n$ with $n>1$.
Therefore,  in contrast to Theorem~\ref {intro-main3}
which yields nonconnected moduli spaces of nonnegative sectional
curvature metrics in all dimensions $\ge 22$, 
it is not clear if the components of  
$\mathfrak R_{\sec \ge 0}(M)/\Diff (M)$ 
remain disconnected after such stabilization.

There are  many other interesting results
about the connectedness or disconnectedness 
of moduli spaces of metrics 
satisfying certain geometric bounds,
for which we refer to, e.g., \cite{L,NW,PRT,PT}.

We conclude the introduction 
with a short description of the ideas and outlines of the proofs.

To prove Proposition~\ref{intro-curv} we look at a  $6$-manifold $X^6$ 
which is homotopy equivalent to $S^2\times S^2\times S^2$ 
but has  nontrivial first Pontrjagin class.  
By an easy topological argument, among the $S^1$-bundles over $X^6$ 
there are infinitely many spaces which are homotopy equivalent to
$S^2\times S^2\times S^3$ but have distinct Pontrjagin classes.
All $S^1$-bundles we consider can be represented 
as quotients of a fixed manifold $Q$ by various subtori 
$T^2_i\subset T^3$ where $T^3$ acts freely and isometrically on $Q$. 
This  implies that the induced metrics on $Q/T^2_i$ have uniformly bounded
curvatures and diameters.

To prove Theorem~\ref{intro-main} 
we fix a rank 2 bundle $\xi$ over $S^2\times S^2\times S^2$ 
and look at the sphere bundle $P$ of $\xi\oplus \e^{k-1}$ with $k\ge 3$. 
We then look at various circle bundles $S^1\to M_i\to P$. 
A topological argument shows that with an appropriate choice of $\xi$,  
infinitely many such bundles
have total spaces homotopy equivalent to $S^2\times S^2\times S^3\times S^k$ 
but distinct first Pontrjagin classes and thus are mutually non-homeomorphic.   

We can represent all $M_i$s as $S^3\times S^3\times S^3\times S^k/T^2_i$ 
where $T^2_i\subset T^3$ 
which acts freely and isometrically on $S^3\times S^3\times S^3\times S^k$. 
This easily implies that the $M_i$ satisfy all geometric constraints 
in Theorem~\ref{intro-main}.

To prove Theorem~\ref{intro-main2} 
we put $k=3$, fix a rank 2 bundle $\zeta$ over $P$ and look at the pullbacks
of $\zeta\oplus \e^{l-2}$  to
$M_i$.  By the same reasons as before, 
the total spaces of these  pullbacks have metrics satisfying all 
geometric restrictions of Theorem~\ref{intro-main2} with souls isometric to $M_i$. 
Another topological argument then shows that with an appropriate choice of 
$\zeta$ the total spaces of the pullbacks are  diffeomorphic to
$S^2\times S^2\times S^3\times S^3\times \R^l$ if $l>10$.

To prove Theorem~\ref{intro-main3}  we modify the construction in the proof of Theorem~\ref{intro-main2}
to produce a manifold with infinitely many nondiffeomorphic souls whose normal bundles have nontrivial rational Euler classes. We then show that for such a manifold all elements of a connected component of $\mathfrak R_{\sec \ge 0}(M)/\Diff(M)$ have diffeomorphic souls.

It is our pleasure to thank Igor Belegradek  
for many helpful discussions and 
particularly for his help with the proof of the 
topological part of
Theorem~\ref{intro-main3}.
We would  also like to thank 
Burkhard Wilking and Wolfgang Ziller 
for helpful conversations regarding the preparation of this article.

   \section{Proof of Proposition~\ref{intro-curv}}

\begin{proof}

Let $\a,\b,\ga$ be the standard basis of $H^2(S^2\times S^2\times S^2)$. 
By Lemma~\ref{lem-hom} in the appendix, for some $m>0$
there exists a closed manifold $M^6$ and a  smooth homotopy equivalence
$f\co M\to S^2\times S^2\times S^2$ such that $p_1(M)=f^*(m\beta\wedge\gamma)$.
   
Consider the principal $T^3$ bundle 
$T^3\overset{p}{\to} S^3\times S^3\times S^3\to S^2\times S^2\times S^2$ 
and let $Q=f^*(p)$ be its pullback. 
Choose a Riemannian metric $g$ on $Q$ which is invariant under the $T^3$ action. 
For any subtorus $T^2\subset T^3$, 
the quotient space $Q/T^2$ is naturally 
a principal $S^1$-bundle over $Q/T^3=M^6$. 
Clearly, any  principal $S^1$- bundle over $M^6$ 
with indivisible Euler class can be realized in this way. 
Let us denote the subtorus corresponding to the bundle with Euler class
 $(a,b,c)$ by $T^2_{a,b,c}$. 
Here the Euler class is written with respect to the  natural product basis 
$\a,\b,\ga$ of $ H^2(M^6)\cong H^2(S^2\times S^2\times S^2) $. 
   
By Lemma~\ref{lem-es}, 
 all the quotients $Q/T^2_{a,b,c}$ with the induced submersion metrics 
satisfy $|\sec|\le C, \diam \le D$ for some $C,D>0$.

Also from Lemma A.3, we see that all the spaces $N_{a,b}=Q/T^2_{a,b,0}$ 
 with $(a,b)=1$ are homotopy equivalent to $S^2\times S^2\times S^3$.

Now, for  $\pi\co N_{a,b}\to M$ with $(a,b)=1$ we have
 $\pi^*(\a)=-b\omega, \pi^*(\b)=a\omega, \pi^*(\ga)=\ga$ 
and thus $\pi^*(\b\wedge\ga)=a\omega\wedge\ga$, 
where $\omega\wedge\ga$ is the generator of
$H^4(N_{a,b})$. Therefore, 
$p_1(N_{a,b})=\pi^*(p_1(M))=am\omega\wedge\ga$.  
This means that all manifolds $N_{a,b}$ with distinct $a$ and $(a,b)=1$  
have distinct Pontrjagin classes and thus are mutually non-homeomorphic.

Finally observe that crossing the manifolds $N_{a,b}$ with round spheres 
produces examples satisfying the conclusion of Theorem~\ref{intro-curv} 
in all dimensions $\ge 8$.
   \end{proof}

\section{Proof of Theorem~\ref{intro-main}}

\begin{proof}

Fix $k\ge 3$.

Consider the standard free $T^3$ action on $S^3\times S^3\times S^3$ 
giving rise to the bundle
$T^3\to S^3\times S^3\times S^3\to S^2\times S^2\times S^2$.  
For any subtorus $T^2\subset T^3$, 
the homogeneous space $S^3\times S^3\times S^3/T^2$ is naturally 
a principal $S^1$-bundle over $S^3\times S^3\times S^3/T^3=S^2\times S^2\times 
S^2$. 
Clearly, any  principal $S^1$- bundle over $S^2\times S^2\times S^2$ 
with indivisible Euler class can be realized in this way. 
Let us denote the subtorus corresponding to the bundle with Euler class
 $(a,b,c)$ by $T^2_{a,b,c}$. 
Here the Euler class is written with respect to the  natural product basis 
$\a,\b,\ga$ of $H^2(S^2\times S^2\times S^2)$. 
Let $N_{a,b,c}$ be the corresponding total space and 
$\pi\co N_{a,b,c}\to S^2\times S^2\times S^2$ be the natural projection.

By Lemma~\ref{lem-easy} , $N_{a,b,0}$ is homotopy equivalent to
 $S^2\times S^2\times S^3$ if $(a,b)=1$.
 
Let us fix a representation $\rho\co T^3\to SO(2)$ 
given by the weight $(p,q,r)$. Look at the associated
$\R^2$ bundle $\xi$ over $S^2\times S^2\times S^2$ given by 
$S^3\times S^3\times S^3\times _{T^3}\R^2$.
Its Euler class is $(p,q,r)$.  
Let $\eta=\xi\oplus \e^{k-1}$ and let $\eta^S$ be 
the corresponding {\it sphere} bundle 
$S^k\to P \overset{\eta^S}{\to} S^2\times S^2\times S^2$ (here  and in what follows $\e^m$ denotes a trivial $\R^m$-bundle).

Next look at the pullback of $\eta$ to $N_{a,b,0}$. 
It can be written as $S^3\times S^3\times S^3\times _{T^2_{a,b,0}}\R^2\times 
\R^{k-1}$.  
We will denote this bundle by $\eta_{a,b}$.  
Let $S^k\to M_{a,b}\overset{\eta^S_{a,b}}{\to}N_{a,b,0}$ 
be the corresponding {\it sphere} bundle over $N_{a,b,0}$.

We claim that by choosing  an appropriate $\rho\co T^3\to SO(2)$ 
and by varying $a,b$ the manifolds $M_{a,b}$ provide examples 
satisfying the conclusion of Theorem~\ref{intro-main}.

Let us first check the geometric conditions.

Observe that  we can write $\eta$ and $\eta^S$ as 
$S^3\times S^3\times S^3\times_{T^3} \R^{k+1}$,
$S^3\times S^3\times S^3\times_{T^3}  S^{k}$ respectively.  
Here $T^3$ acts on $S^3\times S^3\times S^3$ by the canonical homogeneous action 
and on $\R^{k+1}$ and $S^k$  via $\rho$ followed by the canonical inclusion 
$SO(2)\hookrightarrow SO(k+1)$.

Hence, $M_{a,b}=S^3\times S^3\times S^3\times_{T^2_{a,b,0}}  
S^{k}$.

Therefore, by Lemma~\ref{lem-es}, when equipped with the induced quotient metrics, 
all total spaces have uniform curvature bounds $0\le \sec \le C$ 
for some $C>0$ and $\diam\le D$ for some $D>0$.

Next let us check that $\Ric (M_{a,b})>0$.  
Obviously, $\Ric (M_{a,b})\ge 0$. 
Suppose there exists $x_0\in T_pM_{a,b}$ such that $\Ric(x_0)=0$.  
Let $\tilde{x}_0$ be its horizontal lift to 
$T_{\tilde{p}}(S^3\times S^3\times S^3\times S^k)$.  
By O'Neill's formula this means that $\sec(\tilde{x}_0, x)=0$ 
for any horizontal vector 
$x\in T_{\tilde{p}}(S^3\times S^3\times S^3\times S^k)$.
Let $\h$ and $\m$, respectively, denote  the horizontal and the 
vertical tangent space at $\tilde{p}$.

Then $\tilde{x}_0$ contains a nontrivial component tangent to some  sphere factor.
By construction, the projection of $\m$ 
to the tangent space to that sphere is at most one dimensional.
Therefore we can find a vector $x$ tangent to that spherical 
factor and perpendicular to both $\m$ and $\tilde{x}_0$. 
Then $\sec(\tilde{x}_0,x)>0$ which is a contradiction.

 To finish the proof of Theorem~\ref{intro-main}, it remains to show
that among the spaces $M_{a,b}$ there are infinitely 
many homotopy equivalent but mutually non-homeomorphic ones.

 First we claim that there exists an integer $m$ such that  
for any $\rho$ with weight
 $(mp,mq,mr)$,  {\it all}  spaces $M_{a,b}$ are homotopy equivalent 
to $S^2\times S^2\times S^3\times S^k$ if $a$ and $b$ are reletevely prime.
 
 Look at the sphere bundle  $S^k\to P \overset{\eta}{\to} S^2\times S^2\times S^2$.
 
  Up to fiberwise homotopy equivalences such bundles are classified 
by the homotopy classes of maps in  $[S^2\times S^2\times S^2, B\Aut (S^k)]$. 
Here $B\Aut (S^k)$ is the classifying space for $\Aut (S^k)$ 
which is the identity component of the monoid of self-homotopy equivalences of 
$S^k$.
  
  Moreover, in our case, by construction, the classifying map 
into $B\Aut (S^k)$ corresponding to the bundle $\eta$ 
factors through $B\Aut_0 (S^3)$ 
where $\Aut_0(S^k)$ is the subset of $\Aut (S^k)$ fixing a  base point.

 It is a well known fact that if $k$ is odd, then
$\pi_i(\Aut_0(S^{k}))$ is finite for any $i$.
Indeed,  it is easy to see that $\Aut_0(S^k)$ is the identity component 
of $\Omega^k(S^k)$, and therefore, for any $i>0$, 
$\pi_i(\Aut_0(S^k))\cong \pi_{k+i}(S^k)$,  which is always finite   if $k$ is odd.

 A standard obstruction theory argument now implies 
that $[S^2\times S^2\times S^2, B\Aut_0 (S^3)]$ is finite.  
 
 {\bf Claim.}  There is an $m>0$ such that if $e(\xi)$ is divisible by $m$,
 then the classifyng map $f_\eta\co S^2\times S^2\times S^2\to B\Aut_0 (S^3)$ 
is homotopic to a point.
 
 For any $m>0$ let $g_m\co S^2\to S^2$ be a map of degree $m$.  
 
 Let $F_m=g_m\times g_m\times g_m\co S^2\times S^2\times S^2 \to S^2\times 
S^2\times S^2$. Clearly,
 $e(F_m^*\xi)=m e(\xi)$ for any rank $2$ bundle $\xi$ over $S^2\times S^2\times 
S^2$. 
Hence if $f_\eta\co S^2\times S^2\times S^2 \to B\Aut_0 (S^3) $ 
is the classifying map for the $S^k$ bundle coming from $\xi$,
 then $f_\eta\circ F_m$ is the classifying map 
for the $S^k$ bundle coming from the rank 2 bundle with Euler class equal to $m 
e(\xi)$.
 
 The claim now follows from a standard obstruction theory argument. 
 Let us give a brief sketch. 
If $f_1,f_2\co S^2\times S^2\times S^2\to B\Aut_0 (S^3)$ 
are two maps which are homotopic on the $(i-1)$-skeleton, 
the obstruction to extending this homotopy to the $i$-skeleton lies
 in $\G_i=H^i(S^2\times S^2\times S^2, \pi_i(B\Aut_0 (S^3))$ 
which is finite by what has been said above. 
Let $m_i=|\G_i|$.  
By naturality, 
the obstruction corresponding to the maps $f_1\circ F_{m_i}$, 
 $f_2\circ F_{m_i}$ is zero. 
Repeating this process finitely many times, 
we see that for $m=m_1\cdot\ldots \cdot m_6$, 
and any $f_1, f_2$, the maps $f_1\circ F_{m}$, 
 $f_2\circ F_{m}$ are homotopic. 
By taking $f_1$ to be a constant map 
we see that  for any $f\co S^2\times S^2\times S^2\to B\Aut_0 (S^3)$,
the map $f\circ F_{m}$ is homotopic to a constant. This proves our claim.
 
 From now on we will assume that $e(\xi)$ is divisible by $m$ and hence the bundle 
  $S^k\to P \overset{\eta^S}{\to} S^2\times S^2\times S^2$  is 
fiberwise  homotopically trivial.
  
  Of, course the same is true for any pullback of this bundle and hence 
$\eta^S_{a,b}$ is fiberwise  homotopically trivial for any $a,b$. 
Thus its total space $M_{a,b}$ is homotopy equivalent to $N_{a,b,0}\times S^k$ 
which, by Lemma~\ref{lem-easy}, 
is homotopy equivalent to $S^2\times S^2\times S^3\times S^k$ if $(a,b)=1$.

Let us finally show that for appropriately chosen $p,q,r$, 
infinitely many of the spaces $M_{a,b}$  have distinct 
Pontrjagin classes and thus are mutually not diffeomorphic.

Consider the bundle $\pi\co S^1\to Q\to S^2\times S^2$ with
 Euler class $(a,b)$ with respect to the 
canonical generators $\a,\b$ of $H^2(S^2\times S^2)$. 
Let $\omega$ be the generator of $H^2 (Q)$. 
Then  we see from the Gysin sequence that $\pi^*(\a)=-b\omega, \pi^*(\b)=a\omega$.   

Now look at the bundle $\pi\co M_{a,b}\to P$.
Let $\omega,\gamma$ be the natural basis of $H^2(M_{a,b})$. 
Then by the above we have that $\pi^*(\a)=-b\omega,\pi^*(\b)=a\omega, 
\pi^*(\ga)=\ga$.
(We purposefully slightly abuse notations 
by denoting by $\ga$ elements of both $H^2(P)$ and $H^2(M_{a,b})$).

We compute 
\[\begin{split}p_1(\zeta_{a,b})=p_1(\xi_{a,b})
=e(\xi_{a,b})\cup e(\xi_{a,b})=\pi^*(p\a+q\b+r\ga)\cup \pi^*(p\a+q\b+r\ga)
=\\ =((-pb+qa)\omega+r\ga)\cup ((-pb+qa)\omega+r\ga)
=2(-pb+qa)r\omega\wedge\gamma\end{split}\]

 Notice that $\omega\wedge\gamma$ is the generator of $H^4(M_{a,b})\cong Z$. 
 
 From the bundle $S^k\to M_{a,b}\overset{\eta^S_{a,b}}{\to} S^2\times S^2\times 
S^3$,
  using the Whitney formula 
we see that 
$p_1(M_{a,b})=p_1(\eta_{a,b})+\eta^{S*}_{a,b}(p_1(S^2\times S^2\times S^3))=
 2(-pb+qa)r\omega\wedge\gamma$.
 
 Clearly, for fixed $p,q,r$, 
infinitely many of these spaces have distinct $p_1$. 
For example if $p=0, q=mq_1\ne 0,r=mr_1\ne 0$,  
the spaces $M_{a,b}$ with distinct $a$ will work, as in this case 
$p_1(M_{a,b})=2qar$.
 
 By the above, 
all of these spaces are homotopy equivalent  to $S^2\times S^2\times S^3\times 
S^k$
 and hence they satisfy all conclusions of Theorem~\ref{intro-main}.

\end{proof}

\begin{rmk}
   Observe that unlike the examples constructed in ~\cite{Bel}, 
the manifolds constructed in the proof of 
Theorem~\ref{intro-main} have infinite $\pi_2$. 
This is actually necessary by the $\pi_2$-finiteness theorem~\cite{PT}. 

Also, all our examples are constructed as quotients of a fixed manifold $E$  (in our case  $E=S^3\times S^3\times S^3\times S^k$)  by free torus actions. This is also necessary by~\cite[Corollary 0.2]{PT}. 
\end{rmk}

\section{Proof of Theorem~\ref{intro-main2}}\label{sect:souls}
 We will use the notation 
and constructions employed in the proof of Theorem~\ref{intro-main} .

We will make use of the following fact from algebraic topology 
which follows  from a combination of results of Haefliger and Siebenmann 
(~\cite{Hae,Sie}).

{\bf Fact:} 
Let $\R^l\to E_i\to M^n_i, (i=1,2)$ be two vector bundles over smooth closed 
manifolds. 
Suppose $f\co E_1\to E_2$ is a tangential homotopy equivalence and $l\ge 3, l>n$.
   
Then $f$ is homotopic to a diffeomorphism (cf.~\cite{Bel} for details).

Let in the proof of Theorem~\ref{intro-main} now $k=3$.
   
Recall that from the construction of  
$P$ as $S^3\times S^3\times S^3\times_{T^3} S^3$, 
we see that the map $\Z^3\cong\pi_2(P)\to\pi_1(T^3)\cong \Z^3$ is an isomorphism.

Consider a representation $\phi\co T^3\to SO(2)$. 
It gives rise to a rank $2$ bundle $\zeta$ over $P$. 
Its total space can be written as  
$E_\zeta= S^3\times S^3\times S^3\times S^3\times_{T^3} \R^2$, 
where $T^3$ acts on $S^3\times S^3\times S^3\times S^3$ 
by the  action  described above and on $\R^2$ by $\phi$.
   By above, by choosing appropriate $\phi$ we can realize 
in this way {\it any} rank 2 bundle over $P$ with  Euler class $(x,y,z)$  
with respect to the basis $\a,\b,\ga$ of $H^2(P)\cong H^2(S^2\times S^2\times 
S^2)$.
   
   Recall  that by the proof of Theorem~\ref{intro-main},
we can chose $\xi$ so that  $e(\xi)=(0, q, r)$.

Let us choose $\zeta$ so that $e(\zeta)=(\eta^{S})^*(0,q, -r)$ 
where we recall that $\eta^S$ is the sphere bundle 
$S^3\to P\to S^2\times S^2\times S^2$.    
   
Let $\zeta_{a,b}$ be the pullback of $\zeta$   to $M_{a,b}$ 
via the natural projection $\pi\co M_{a,b}\to P$.
 
We will show that infinitely many of the stabilized bundles 
 $\tilde{\zeta}_{a,b}=\zeta_{a,b}\oplus \e^{l-2}$ 
satisfy the statement of Theorem~\ref{intro-main2} if $l>10$.
 
 Let us first check the geometric conditions. 
 The total space of $\tilde{\zeta}_{a,b}$ can be written as
 $ E({\tilde{\zeta}_{a,b}})=S^3\times S^3\times S^3\times S^3\times_{T^2_{a,b}} 
\R^l$. 
Here $T^2_{a,b}\subset T^3$ is the subtorus which corresponds 
to the bundle $S^1\to M_{a,b}\to P$. 
Therefore, by Lemma~\ref{lem-es},  we have uniform curvature bounds.
Note that while the manifolds in question are not compact, 
it is easy to see that curvature remains uniformly bounded at infinity, 
so that Lemma~\ref{lem-es} still applies.
Alternatively, rather than taking $\R^l$ with a flat metric, 
we can take it with a rotationally symmetric nonnegatively 
curved metric {\it isometric } to $S^{l-1}\times \R_+$ at infinity. 
Then the uniform curvature bounds follow directly from Lemma~\ref{lem-es}.
  
 Of course, the soul of $E(\tilde{\zeta}_{a,b})$ is isometric to $M_{a,b}$,
 and thus all the souls have bounded diameter and are not homeomorphic for 
different $a$.
   
Next we will show  that infinitely many of the bundles $\tilde{\zeta}_{a,b}$
have diffeomorphic total spaces.
    
First, by the same computation as in the proof of Theorem~\ref{intro-main}, we 
find
$p_1(\tilde{\zeta}_{a,b})=-2qar$ and hence $p_1(E((\tilde{\zeta}_{a,b}))
=p_1(\tilde{\zeta}_{a,b})+p_1(M_{a,b})=-2qar+2qar=0$.
    
Since Pontrjagin classes determine a bundle up to finite ambiguity, 
infinitely many of the spaces $E({\tilde{\eta}_{a,1}})$ 
are  tangentially homotopy equivalent and hence diffeomorphic.

This observation  is already sufficient to produce examples of manifolds with infinitely many nonnegatively pinched metrics whose souls have bounded diameter and are mutually non-homeomorphic. Unfortunately, it does not give us the precise diffeomorphism type of these manifolds.

However,  with a little more work we can show that, in fact, 
{\it all} manifolds 
$E({\tilde{\zeta}_{a,b}})$ with $(a,b)=1$ 
are diffeomorphic to $S^2\times S^2\times S^3\times S^3\times \R^l$.

Look at the following commutative diagram:

\[
\xymatrix{
S^3\ar[r] & M_{a,b}\ar[r]^{\eta^S_{ab}}\ar[d]_\pi &N_{a,b,0}\ar[d]_\pi\\
S^3\ar[r] & P\ar[r]^-{\eta^S}& S^2\times S^2\times S^2
}
\]

First notice that the bundle $\zeta$ is the pullback via $\eta^S$ 
of the bundle $\hat{\zeta}$ over $S^2\times S^2\times S^2$ with "the same" Euler 
class. 
Bundle $\hat{\zeta}$ can be written as 
$S^3\times S^3\times S^3\times_{T^3}\R^2$ 
where $T^3$ acts on $\R^2$ by the representation $\phi$.

Similarly, $\zeta_{a,b}=\eta_{a,b}^{S*}(\hat{\zeta}_{a,b})$.

Next observe that 
$TE(\tilde{\zeta}_{a,b})|_{M_{ab}}=TM_{a,b}\oplus \zeta_{a,b}\oplus \e^{l-2}=
TM_{a,b}\oplus \e^1\oplus \zeta_{a,b}\oplus \e^{l-3}=
\eta^{S*}_{a,b}TN_{a,b,0}\oplus 
\eta_{a,b}^{S*}(\xi_{a,b}\oplus \e^2) 
\oplus \zeta_{a,b}\oplus \e^{l-3} =
\eta_{a,b}^{S*}(TN_{a,b,0}\oplus \xi_{a,b}\oplus \hat{\zeta}_{a,b}\oplus \e^{l-
1})$.

Since $N_{a,b,0}$ is the total space of an $S^1$ bundle over 
$S^2\times S^2\times S^2$,  it immediately follows that 
$TN_{a,b,0}\oplus \e^{l-1}=e^{l+6}$. (Alternatively, 
this is also clear since $N_{a,b,0}$ is 
diffeomorphic to $S^2\times S^2\times S^3$ by~\cite{Ba}).

  Thus $TE(\tilde{\zeta}_{a,b})|_{M_{ab}}=\eta_{a,b}^{S*}(\xi_{a,b}
\oplus \hat{\zeta}_{a,b})\oplus \e^{l+6}=
\eta_{a,b}^{S*}\pi^*(\xi\oplus\hat{\zeta})\oplus \e^{l+6}=
\eta_{a,b}^{S*}\pi^*(\e^4)\oplus \e^{l+6}=\e^{10+l}$. 
Here the next to last equality holds 
because  by the choice of $\zeta$ we have that 
$e(\hat{\zeta})=(0,q,-r)$ and $e(\xi)=(0,q,r)$ 
so that $\xi\oplus\hat{\zeta}=\e^4$ by Lemma~\ref{lem2} below.

Thus $E(\tilde{\zeta}_{a,b})$ is tangentially equivalent 
and hence diffeomorphic to $S^2\times S^2\times S^3\times S^3\times \R^l$.

\hfill $\qed$
   
\begin{rmk}
Using the same procedure  as in the proof of Theorem~\ref{intro-main2}, 
 we can also construct manifolds with nontrivial $p_1$ 
which admit  infinitely many nonnegatively pinched metrics with non-homeomorphic souls.
\end{rmk}

\section{Proof of Theorem~\ref{intro-main3}}
We will use the same notations as in the proofs of Theorem~\ref{intro-main} and Theorem~\ref{intro-main2}.
Let us first construct the Riemannian manifolds in question. The construction is very similar to the one used in the proof of Theorem~\ref{intro-main2}, therefore we will skip some details.

Let $m$ be as in the proof of Theorem~\ref{intro-main}.
Let us fix positive integers $n$ and $k$ and look at a 
rank 2 vector bundle $\xi$ over $S^2\times S^2\times S^2$ with 
Euler class $m(0,1,k)$. 
Let $S^3\to P\to S^2\times S^2\times S^2$ 
be the sphere bundle in $\xi\oplus \e^2$. 
Look at the rank 2 bundles $\zeta_1, \zeta_2$ over $P$ with  
Euler classes $m\pi^*(0,1,-k), m\pi^*(n+1, n, 0)$, respectively.
Now look at the $S^1$ bundle $\pi_{ab}\co N_{a,b}\to P$ over $P$ 
with  Euler class $\pi^*(a,b,0)$, where $(a,b)=1$ and pull back
 $\zeta=\zeta_1\oplus \zeta_2 $ to $N_{a,b}$. 

By the proof of Theorem ~\ref{intro-main}, $N_{a,b}$ is homotopy equivalent to $S^2\times S^2\times S^3\times S^3$ for any 
pair of integers $a,b$ with $(a,b)=1$.

As before we also see that the total space 
of the bundle $\pi_{a,b}^*(\zeta)$ admits a complete metric with
$0\le \sec\le C$ and  the soul isometric to $N_{a,b}$ with $\diam (N_{a,b})\le D$ where $C,D$ are independent of $a,b$.

A computation similar to the one  
in  the proof of Theorem~\ref{intro-main2}
shows that for the first Pontrjagin and Euler classes of the bundles in
question we have
\[
p_1(N_{a,b})=2m^2 ak\omega\wedge\gamma , \qquad
p_1(\pi_{a,b}^*(\zeta))=-2m^2 ak\omega\wedge\gamma ,
\]
and 

\[
e(\pi_{a,b}^*(\zeta))=m^2(-b(n+1)+an)k\omega\wedge\gamma .
\]

Set $a=1+r(n+1), b=1+rn$ where $r\in \N$ and let $N_r=N_{1+r(n+1), 1+rn}$.

Then 

\begin{equation}\label{eq:1}
p_1(N_{r})=2m^2 ak\omega\wedge\gamma , \qquad
p_1(\pi_{r}^*(\zeta))= -2m^2a k\omega\wedge\gamma ,
\end{equation}
and 

\begin{equation}\label{eq:2}
e(\pi_{r}^*(\zeta))=-m^2k\omega\gamma .
\end{equation}

This means that the manifolds $N_r$ have distinct Pontrjagin classes 
and hence are mutually non-homeomorphic.

Let $E_r$ be the total space of the bundle $\pi_{r}^*(\zeta)$. 
From the above we see that 

\begin{equation}\label{eq:3}
p_1(E_r)=2m^2 ak\omega\wedge\gamma-2m^2a k\omega\wedge\gamma=0 .
\end{equation}

Look at the spaces $X_r=E_r\times TS^4$ with the product metric where
we take the natural nonnegatively curved metric on $TS^4$ given by the submersion metric on $TS^4=SO(5)\times_{SO(4)}\R^4$.

We claim that  
the spaces  $X_r$ fall into 
finitely many diffeomorphism classes.

Indeed, let $f_r\co N_{1}\to E_{r}$ be the homotopy equivalence 
given by the homotopy equivalence of the souls 
followed by the embedding of the soul into $E_{r}$.  
Note that  $H^4(N_r)\cong H^4(S^2\times S^2\times S^3\times S^3)\cong \Z$. By (\ref{eq:2}), by possibly composing $f_r$ with an 
orientation reversing self homotopy equivalence of $N_{1}$,
 we can assume that $f_{r}^*(e(\pi_{r}^*(\zeta)))=e(\pi_{1}^*(\zeta))$.

Observe that $X_{r}$ is the total space of a rank $4+4=8$ vector bundle over $N_{r}\times S^4$ and $\dim N_r\times S^4= 10+4=14$. Since $f_r\times \Id_{S^4}$ is a homotopy equivalence and $3\cdot 8>14+2$, 
we are in the metastable range and by Haefliger's  
Embedding Theorem~\cite{Hae},  $f_r\times \Id_{S^4}$ 
is homotopic to an embedding $g_r$. Since the codimension of 
$N_r\times S^4$ in $X_r$ is $=8>3$, by~\cite{Sie}, 
$X_{r}$ is diffeomorphic to the total space of the 
normal bundle $\nu_{g_r}$. 

From (\ref{eq:1}) and (\ref{eq:3}), using the 
Whitney formula we see that all $\nu_{g_r}$ 
have the same Pontrjagin classes.
From (\ref{eq:2}) we see also 
that all $\nu_{g_r}$ have the same {\it nontrivial} 
Euler classes equal to $-2m^2k \omega\wedge \gamma\wedge [d_{vol}(S^4)]$. That is because the rational  Euler class of $\nu_{g_r}$ is a 
homotopy invariant of $g_r$ which can be defined homologically 
by the formula
$<e(\nu_{g_r}), x>=g_{r*}[X_1]\cdot g_{r*}(x)$  for any  $x\in H_8(X_1)$ where $\cdot$ is the algebraic intersection number. Thus
 $e(\nu_{g_r})=f_{r}^*(e(\pi_{r}^*(\zeta))\cup e(S^4))=
-m^2k\omega\gamma\wedge 2 [d_{vol}(S^4)]$ by (\ref{eq:2}) 
and the fact that $\chi(S^4)=2$. 
Here we disregard the difference 
between rational and integer coefficients 
since all involved  cohomology groups  are torsion free.
See also ~\cite{BK} for a   more detailed discussion of invariants of maps.

Thus we see that all the bundles  $\nu_{g_r}$ have the same Euler and Pontrjagin classes. Since Euler and Pontrjagin classes determine a bundle up to a finite ambiguity, the bundles $\nu_{g_r}$ fall into finitely many isomorphism classes. Hence the total spaces of $\nu_{g_r}$ fall into finitely many diffeomorphism classes. By the 
above, the total space of $\nu_{g_r}$ is diffeomorphic to $X_{r}$ 
and hence all manifolds $X_{r}$  
also  fall into finitely many diffeomorphism classes.

Thus after passing to a subsequence, we can assume that all $X_{r_i}$ are diffeomorphic to $M=X_{r_1}$.

We claim that $M$ satisfies the conclusion of the Theorem.
Observe that $X_{r}$ carries by construction 
a natural metric of $0\le sec\le C$ with soul isometric to $N_r\times S^4$ of $\diam \le D$. Hence  all the souls have distinct Pontrjagin classes by (\ref{eq:1}) and thus are mutually not homeomorphic.

Since for any nonnegatively curved metric $g$ and any 
self-diffeomorphism of the underlying manifold $\phi$, the souls of $g$ and $\phi^*(g)$ are diffeomorphic, the statement of Theorem~\ref{intro-main3} will follow from the following
\begin{lem}\label{l:souls}
Let $(M,g_t), t\in[0,1]$ be a continuous family of nonnegatively curved metrics such that  
the normal bundle to the soul of $(M, g_0)$ has 
nontrivial rational Euler class.

Then all the souls of $(M,g_t)$ are diffeomorphic.
\end{lem}
\begin{proof}
Let $S_t$ be the soul of $(M,g_t)$. We claim that the family $(S_t, g_t|_{S_t})$ is continuous in Gromov-Hausdorff topology. Observe that since $S_t\hookrightarrow M$ is a homotopy equivalence, by the same argument as above the rational Euler class of $\nu_{S_t}$ is nonzero for any $t$. Therefore it's enough to show that $S_t\overset{G-H}{\to} S_0$ as $t\to 0$.

Let $\pi_t\co M\to S_t$ be the Sharafutdinov retraction with respect to $g_t$.

Let $d_t$ be the inner metric on $M$ induced by $g_t$.
Since $g_t\to g_0$ uniformly on compact sets we clearly have that for any $x,y\in S_0$, $d_t(x,y)\le d_0(x,y)+\e_t$
 where $\e_t\to 0$ as $t\to 0$.   Since $\pi_t$ is distance nonincreasing we see that $d_t(\pi_t(x),\pi_t(y))\le d_0(x,y)+\e_t$ for any $x,y\in S_0$.  Since $\pi_t\co S_0\to S_t$ is a homotopy equivalence, it must be onto and hence
 $\diam S_t\le \diam S_0+\e_t$.
 
From 
the assumption on the Euler class we see that $S_t\cap S_0\ne \emptyset$ for any $t$ and  since,  by the above, 
all $S_t$  have uniformly bounded diameters, they all must lie in some fixed closed ball $\bar{B}(p,D)$ where the ball is taken with respect to $d_0$.
 Again using that  $g_t$ converges to $g_0$ uniformly on compact sets we have that  $d_0(x,y)\le d_t(x,y)+\e_t$ for any $x,y\in S_t$. Hence  $d_0(\pi_0(x),\pi_0(y))\le d_t(x,y)+\e_t$ for any $x,y\in S_t$. Combining this with above we finally get that
 
 \[
 d_0(\pi_0(\pi_t(x)), \pi_0(\pi_t(x)))\le d_0(x,y) +2\e_t \quad \text{ for any } x,y \in S_0 .
 \]
 
 By Lemma~\ref{lem:haus} this implies that
for some $\tilde{\e} (t)\underset{t\to 0}{\to} 0$
 
  \[
 d_0(x,y) -2\tilde{\e}_t \le d_0(\pi_0(\pi_t(x)), \pi_0(\pi_t(x)))\le d_0(x,y) +2\e_t \quad \text{ for any } x,y \in S_0 .
 \]

 Hence $\pi_0\circ \pi_t\co S_0\to S_0$ is a $\max(\e_t,\tilde{\e}_t)$-Hausdorff approximation and the same is true for
 $\pi_0\co (S_t,d_t)\to (S_0,d_0)$ which proves that $S_t\overset{G-H}{\to} S_0$ as $t\to 0$.
 
 Since $S_t$ is a smooth manifold for any $t$ and $\dim S_t=\dim S_0$, by Yamaguchi's Stability theorem~\cite{Yam}  this implies that
 $S_t$ is diffeomorphic to $S_0$ for all small $t$.
\end{proof}

As observed before, Lemma~\ref{l:souls} implies that all elements of a  connected component of $\frak R_{\sec\ge 0}(M)/\Diff (M)$ have 
diffeomorphic souls. This immediately implies 
the statement of Theorem~\ref{intro-main3}. \hfill $\qed$
      
      \begin{rmk}
      We suspect that Lemma~\ref{l:souls} is true without any assumptions on the rational Euler class. If this holds true,
then the examples constructed in the proof of Theorem~\ref{intro-main2} would directly yield Theorem~\ref{intro-main3}.
      \end{rmk}
\appendix \section{}

   We will need the following  lemma 
which is an easy consequence of some well-known topological results:
   
   \begin{lem}\label{lem-hom}
   There exists an integer $m$ such that for any element 
$p\in H^4(S^2\times S^2\times S^2)$ 
there exists a  closed smooth manifold $M^6$ 
and a homotopy equivalence $f\co M\to S^2\times S^2\times S^2$ 
such that $f^*(p)=mp_1(M)$.
   \end{lem}
   \begin{proof}
   
   By the Browder-Novikov Surgery Theorem~\cite[Thm II.3.1,Cor II.4.2]{Bro},  
given a vector bundle $\xi$ over a simply connected manifold $X^6$, 
there exists a  manifold $M^6$ and a homotopy equivalence $M^6\to X$,
such that $f^*(\xi)$ is isomorphic to the stable normal bundle of $M$ if and only 
if
   the stable spherical fibration coming from $\xi$ is isomorphic 
to the Spivak normal spherical fibration
   $\nu(X)$.
   
   If $X=S^2\times S^2\times S^2$, 
we obviously have that $\nu(X)$ is trivial. 
Recall that stable spherical fibrations are classified by the homotopy 
classes of maps into the classifying space $BG$ and that all 
homotopy groups of $BG$ are finite.  The same obstruction theory argument 
as in the proof of Theorem~\ref{intro-main} shows 
that there exists an $m_1$ such that for any 
$f\co S^2\times S^2\times S^2\to BG$, the map $f\circ F_{m_1}$ is homotopic to a 
point.  Recall here that   $F_m=g_m\times g_m\times g_m\co S^2\times S^2\times S^2 \to S^2\times 
S^2\times S^2$ where $g_m\co S^2\to S^2$ has  $\deg g_m=m$.
   
  Next observe that by looking at Whitney sums of rank 2 bundles 
we can realize any {\it even}  element of 
  $H^4(S^2\times S^2\times S^2)$ as the first Pontrajagin class of a vector bundle.

   Combining these two facts, we obtain the desired claim with $m=2m_1$.
   \end{proof}

   The geometric part of the proof of Theorem~\ref{intro-main}
is based on the following 
lemma, which is originally due to Eschenburg~\cite[Prop 22]{E1}.
For convenience of the reader, we include a short outline of its proof.

\begin{lem}\label{lem-es}
Let $(M,g)$ be a closed Riemannian manifold 
on which a $k$-dimensional torus $T^k$ acts freely and isometrically.
Then there exist $C,D>0$ such that for any subtorus $T^m\subset T^k$ 
the quotient manifold $M/T^m$, when equipped with the induced quotient metric, 
satisfies
\[|\sec(M/T^m)|\le C \,  \text{ and } \, \diam (M/T^m)\le D.
\]

\end{lem}

\begin{proof}
The uniform diameter bound is obvious,
and we need to find a uniform bound 
on the O'Neill term in the Gray-O'Neill curvature formula
for Riemannian submersions. 
As the formula is local, it makes sense to look at the local quotients of $M$ 
by $\R^m\subset \R^k$ where $\R^k$ is the universal cover of $T^k$.  
The compactness of the Grassmannian
of $m$-planes in $\R^k$ now implies the result.

\end{proof}

\begin{lem}\label{lem-easy}
Let $S^1\to P\to S^2\times S^2$ be  a principal $S^1$ bundle such that $P$ 
is simply connected.
Then $P$ is homotopy equivalent to $S^2\times S^3$.
\end{lem}

In fact---though we will not need this fact in this paper---by a theorem of 
Barden~\cite{Ba}, $P$ is 
{\it diffeomorphic} to $S^2\times S^3$.

\begin{proof}

It is easy to see that $H_2(P)\cong H_3(P)\cong \Z$.

We can write $P$ as the homogeneous space 
$S^3\times S^3/S^1$ for some $S^1\subset S^3\times S^3$.
>From the Gysin sequence 
it is easy to see that the map $H_3(S^3\times S^3)\to H_3(P)$  is onto. 
Since,  by the Hurewicz theorem, $\pi_3(S^3\times S^3)\cong H_3(S^3\times S^3)$
and since, by the long exact homotopy sequence, 
$\pi_3(P)\cong \pi_3(S^3\times S^3)$ 
we see that 
$\pi_3(P)\to H_3(P)$ is also surjective. 
Let $f\co S^2\to P$ be a map representing a generator of $\pi_2(P)\cong 
H_2(P)\cong \Z$. 
Let $g\co S^3\to P$ be a map representing a generator of  $H_3(P)\cong Z$. 
Let $\hat{g}\co S^3\to S^3\times S^3$ be a lift of $g$. 
Such a lift exists by the previous discussion.  Consider the map
$F\co S^2\times S^3\to P$ given by 
$F(x,y)=\hat{g}(y)\cdot f(x)$, 
where $x\in S^2, y\in S^3$ and where the $\cdot$ 
represents the homogeneous space action of $S^3\times S^3$ on $P$.
It is straightforward to check that 
$F$ induces an isomorphism on homology and thus is a homotopy  equivalence.

\end{proof}

\begin{lem}\label{lem2}
Let $\xi_1,\xi_2$ be rank 2 bundles over 
$S^2\times S^2$ such that $e(\xi_1)=(q,r)$ and $e(\xi_2)=(q,-r)$ 
with respect to the canonical basis of $H^2(S^2\times S^2)$. 

Then $\xi_1\oplus\xi_2$ is trivial.

\end{lem}

\begin{proof}
A direct computation shows 
that $w_2(\xi_1\oplus\xi_2)=p_1(\xi_1\oplus\xi_2)=e(\xi_1\oplus\xi_2)=0$.

Since, $w_2(\xi_1\oplus\xi_2)=0$, 
its classifying map $f$ into $BSO(4)$ factors  as $f=g\circ c$.
where $c\co S^2\times S^2\to S^4$ is the collapsing map of degree 1. 
Thus $\xi_1\oplus\xi_2=c^*(\eta)$ where $\eta$ is a rank 4 bundle over $S^4$. 
Since $c^*$ is an isomorphism on $H^4$, 
$p_1(\eta)=e(\eta)=0$ and hence $\eta$ is trivial.
\end{proof}

\begin{lem}\label{lem:haus}
Let $S$ be a closed Riemannian manifold. There exists a function $\d\co \R_+\to \R_+ $ such that $\d(\e)\to 0$ as $\e\to0$ and 
such that the following holds.
If $f\co S\to S$ is a homotopy equivalence satisfying
\[
d(f(x),f(y))\le d(x,y)+\e \qquad \text { for any } x,y\in S ,
\]
then

\[
d(x,y)-\d(\e)\le d(f(x),f(y))\le d(x,y)+\e \qquad \text { for any } 
x,y\in S .
\]

\end{lem}
\begin{proof}
Suppose Lemma~\ref{lem:haus} is false. Then there exists a 
sequence $f_i\co S\to S$  as well as a sequence $\e_i\to 0$ satisfying
\[
d(f(x),f(y))\le d(x,y)+\e_i \qquad \text { for any } x,y\in S
\]

such that for some $\d>0$ there exist $x_i,y_i\in S$ such that $d(f(x_i),f(y_i))\le d(x_i,y_i)-\d$.
By Arzela-Ascoli and the 
compactness of $S$ we can assume that $f_i$ uniformly converges to $f\co S\to S$ and $x_i\to x_0, y_i\to y_0$. 
Then $f$ is 1-Lipschitz and $d(f(x_0),f(y_0))\le d(x_0,y_0)-\d$.  By uniform convergence $f_i$ is homotopic to $f$ for large $i$. Hence $f$ is onto. A surjective  1-Lipschitz self-map of a closed manifold  has to preserve the volume which easily implies that it must be an isometry. Therefore we must have $d(x_0,y_0)=d(f(x_0),f(y_0))$. This is a contradiction and hence Lemma~\ref{lem:haus} is true.
\end{proof}

\small
\bibliographystyle{alpha}

\end{document}